\documentclass[11pt]{article}

\usepackage{amsmath,latexsym,epsfig}
\usepackage{amsfonts,amssymb,amscd}

\newcommand{\qed}{{\large $\Box$}}

\newcommand{\grgen}{{\xi}}

\title{\bf  A geometric parametrization for the virtual Euler characteristics of the moduli spaces
of real and complex algebraic 
curves
}
\author{I.P.Goulden\thanks{Department of Combinatorics
and Optimization, University of Waterloo, Waterloo, Ontario, Canada, N2L~3G1}, 
J.L.Harer\thanks{Department of Mathematics, Duke University, Durham, NC~27708-0320, U.S.A.},  
and 
D.M.Jackson\thanks{Department of Combinatorics
and Optimization, University of Waterloo, Waterloo, Ontario, Canada, N2L~3G1}}

\date{ January 11, 1999}

\begin{document}
\maketitle

\newtheorem{theorem}{Theorem}[section]
\newtheorem{proposition}[theorem]{Proposition}
\newtheorem{definition}[theorem]{Definition}
\newtheorem{axiom}[theorem]{Axiom}
\newtheorem{lemma}[theorem]{Lemma}
\newtheorem{corollary}[theorem]{Corollary}
\newtheorem{remark}[theorem]{Remark}
\newtheorem{example}[theorem]{Example}
\newtheorem{conjecture}[theorem]{Conjecture}

\def\cM{{\cal{M}}}
\def\cP{{\cal{P}}}
\def\cT{{\cal{T}}}
\def\cV{{\cal{V}}}
\def\cW{{\cal{W}}}

\def\bfi{{\rm\bf i}}
\def\bfu{{\rm\bf u}}
\def\bfx{{\rm\bf x}}
\def\bfy{{\rm\bf y}}
\def\bfM{{\rm\bf M}}

\def\bbP{{\mathbb{P}}}
\def\reals{{\mathbb{R}}}
\def\comp{{\mathbb{C}}}
\def\integ{{\mathbb{Z}}}
\def\tsh{{\textstyle{\frac{1}{2}}}}
\def\tsqu{{\textstyle{\frac{1}{4}}}}
\def\tr{{\rm{trace\,}}}
\def\rintn{\int_{\reals^N}\!\!}
\def\symgp{{\mathfrak{S}}}
\def\av#1{\left\langle#1\right\rangle_{\cW_N}}
\def\efM{e^{-\tsqu\tr\bfM^2}}
\def\etM{e^{-\tsh\tr\bfM^2}}
\def\ip#1#2{\left\langle #1,#2\right\rangle}
\def\avv#1{\left\langle#1\right\rangle_{\cV_N}}
\def\Fix{{\sf{Fix}}}
\def\ipt#1#2#3{\left\langle #1,#2\right\rangle_{#3}}
\def\stat{{\sf{stat}}}
\def\ones{{\bf 1}}
\def\rel{\,{\sf rel}\,\,}
\def\mybar{\overline}

\def\proof{{\rm\bf Proof:\quad}}

\tableofcontents

\section{Introduction}\label{intr}
We show that the virtual Euler characteristics  for the
moduli spaces of $s$-pointed algebraic curves of genus $g$ can be 
determined from a polynomial $\xi^{s}_g(\gamma)$ in $\gamma^{-1}$
where the parameter $\gamma$  permits specialization, through $\gamma=1$, to the
complex case treated by Harer and Zagier and, through $\gamma=1/2,$ to the real case.
The latter yields  a new result for the space of real curves of genus $g$ with 
a fixed point free involution, giving the Euler characteristic as 
$(-2)^{s-1}(1-2^{g-1})(g+s-2)! {B_g}/g!$ where $B_g$ is the $gth$ Bernoulli number.

A detailed analysis of the parametrized quantity $\xi^{s}_g(\gamma)$ provides
startling evidence that it has geometric significance. Indeed, it may in fact be the virtual
Euler characteristic of some moduli spaces, as yet unidentified.
We carry out this analysis by establishing a connection with a previous combinatorial
conjecture that the indeterminate
$b=\gamma^{-1}-1$ is associated with a combinatorial invariant for cell decompositions.
The development uses Strebel differentials to triangulate the moduli spaces, and
the identification of $\gamma$ both as a parameter in  a Jack symmetric function
through the irreducible characters of the symmetric group and the hyperoctohedral
group, and as a parameter in a matrix model through generalized Selberg integrals.
Both the matrix models of Hermitian and real symmetric matrices 
and these character theoretic formulations for counting cell 
decompositions are used.

To explain the parametrization
that we have in mind, we begin with a brief introduction  first to the complex case and then
to the real case.


\subsection{The complex and real cases}\label{SSb}

The moduli space $\cM_g^s$ of $s$-pointed complex curves of genus $g$ has
been an important object of study in many branches of mathematics for
a long time.  The global structure of these spaces is only partially
understood; in particular, the cohomology of $\cM_g^s$ is known in just
a few cases.  The primary tool for studying this cohomology has been a
triangulation of $\cM_g^s$ which is derived from 
work of Strebel on quadratic differentials by Harer~\cite{Harer2}.  In particular,
the Euler characteristic of $\cM_g^s$ was computed by Harer and Zagier~\cite{HZ} using
this triangulation and an associated matrix model to count cells.  Since that
time, much work has been done in applying other matrix models
and variants of the Strebel triangulation to the study of moduli spaces.  This
includes the alternative formulation by Penner~\cite{PE} of the work in ~\cite{HZ}, 
the character-theoretic approach to the same work by Jackson~\cite{JONE}, 
the  exposition by Itzykson and Zuber~\cite{IZ1},
and most recently the beautiful work of Kontsevich~\cite{Kont} in establishing a 
conjecture
of Witten~\cite{Witten} on intersection numbers of cohomology classes of $\mybar{\cM}_g^s$.

The moduli spaces of real algebraic curves have only come into the 
foreground more recently.  
Classically, a real algebraic curve $C$ is the set of simultaneous 
zeros in $\reals \bbP^n$ of a collection of homogeneous polynomials.  By taking the
solution in $\comp \bbP^n$ of the same set of polynomials, we can associate to the
curve a Riemann surface and an antiholomorphic involution on it 
(induced by complex conjugation) whose fixed point set is the real curve $C$.
This suggest {\em defining} a real curve as a pair consisting of a complex 
curve $X$ and an antiholomorphic involution $\tau$ and it is this point of view
that we shall use here.  For a fixed topological type of orientation reversing 
involution $\tau$, there is then a moduli space $\cM_g^s(\tau)$
whose points are isomorphism
classes of pairs of this form.  Sepp\"{a}l\"{a}, Silhol, Buser and others~\cite{BS,S,SS} 
have  made some progress in understanding these spaces, but
up to now little is known about the topology of $\cM_g^s(\tau)$.
In the present work we will apply the technique of the Strebel differential
to obtain a triangulation of $\cM_g^s(\tau)$.  We then count cells using 
enumerative techniques similar to those that have been used to study
$\cM_g^s$. 

It is particularly startling that the complex and real cases can be treated
simultaneously, through a generalized Selberg integral with parameter $\gamma,$
and it is this that leads to the parametrization $\xi^s_g(\gamma)$ of the
Euler characteristic.


\subsection{The geometric parametrization}\label{SScceg}

Central to our approach are the two generating series, $M^O$  for graphs 
embedded in orientable surfaces and  $M$ for graphs embedded in locally orientable
surfaces. 
Collectively, orientable and nonorientable
surfaces are referred to as {\em locally} orientable.  (Of course every
surface is locally orientable, but the terminology is provided to make
it clear that we are considering both orientable and non-orientable surfaces
at once.)
The embedded graphs are combinatorial objects and we will call them {\em maps}.  

$M^O$ and $M$ have been given in integral form as evaluations of expectation 
operators over Hermitian complex matrices and real symmetric matrices. For orientable 
surfaces these are familiar from the work of 't~Hooft~\cite{tH}  and Bessis, Itzykson and 
Zuber~\cite{BIZ}, while the corresponding result for surfaces is given by 
Goulden and Jackson~\cite{GJI}. In the particular case of specialization to the virtual 
Euler characteristics of the moduli spaces, the diagonalized forms of the associated 
expectations have a common generalization $\xi^s_g(\gamma)$ as a generalized Selberg 
integral, whose evaluation is known~\cite{Selb}. The virtual Euler characteristics for  
$\cM_g^s$  and $\cM_g^s(\tau)$
are therefore computed simultaneously by specializing the parameter, $\gamma,$  
that appears in this Selberg integral to $\gamma=1$ and $\gamma=1/2,$ respectively. 

Character theoretic formulations for the generating series $M^O$ and $M$ have been given
by  Jackson and Visentin~\cite{JV} and Goulden and Jackson~\cite{GJ2} 
in terms of Schur functions and zonal polynomials, 
respectively. Both of these functions are specializations of the Jack symmetric
function through its parameter $\gamma.$ We have conjectured~\cite{GJJ}  that,
in the character theoretic setting, $b=\gamma^{-1}-1$ is an indeterminate associated with
a combinatorial invariant of nonorientability for maps, and this identification 
is  therefore inherited in the Selberg integral formulation of this paper.
In the context of moduli spaces, this enables us therefore to
conjecture that the parameter $b$ also has geometric
significance, and that there is an appropriately defined
moduli space involving $b$ that specializes to $\cM_g^s$ and $\cM_g^s(\tau)$ at
$b=0$ and $b=1,$ respectively.

Jack symmetric functions and their dependence on the parameter $\gamma$ 
have been studied quite intensively in the combinatorial community
in the last decade (see, for example, Stanley~\cite{ST}. 
They also appear in connection~\cite{LPS} with the
Calogero-Sutherland model, for example, in statistical mechanics.

\subsection{Organization of the paper}\label{SSop}
The paper is organized as follows. 
Section~\ref{Stsu} gives the reduction of the
calculation of the virtual Euler characteristic $\chi(\cM_g^s(\tau))$ to
a combinatorial sum involving polygonal identifications.
In Section~\ref{Sge} we show that both this combinatorial sum, and the
corresponding sum for the complex case, can be realized as particular coefficients
in the map series $M$ and $M^O$  at transformed arguments.
We then evaluate these by first representing them each as matrix integrals
and then constructing a single generalized Selberg integral $M_\gamma,$ with parameter $\gamma,$ 
that specializes to the diagonalized forms of these matrix integrals through  
choices of $\gamma.$  
This integral is evaluated in Section~\ref{Spmr} through asymptotic properties
of the gamma function to give
the parametrized Euler characteristic,  $\xi_g^s(\gamma),$
explicitly as a polynomial in $\gamma^{-1}.$ 
In Section~\ref{Scrvec} we demonstrate that $M_\gamma$ can be reexpressed
in terms of Jack symmetric functions with parameter $\gamma^{-1}.$
A previous conjecture that gives combinatorial meaning to the parameter $b$
(recall that $b=\gamma^{-1}-1$) in the context of maps
then leads to the conjecture that $b$ has geometric significance
in the context of moduli spaces.
For completeness, we include a brief account of the combinatorial encoding of maps
and the algebras that carry the encoded combinatorial information that
explain the appearance of Schur functions and zonal polynomials
(which are Jack symmetric functions at $\gamma=1$ and $\gamma=1/2,$ respectively).

The Appendix gives a table of the polynomial coefficients in $b$
that arise in the first few terms of $M_\gamma.$
The polynomials are indexed by vertex distribution, the number
of faces and the number of edges as combinatorial data. At $b=0$ and $b=1$ these evaluate
to the numbers of maps, with prescribed combinatorial data, on orientable
and locally orientable surfaces respectively.


\section{The moduli spaces of real curves}\label{Stsu}

\subsection{The Teichm\"{u}ller spaces and mapping class groups of nonorientable surfaces}

By a real algebraic curve we will mean a pair $(X,\sigma)$ where $X$ is a Riemann surface 
and $\sigma$ is an antiholomorphic involution.  Two such curves 
$(X_1,\sigma_1)$ and $(X_2,\sigma_2)$ are isomorphic if there is a biholomorphism 
$\phi: X \rightarrow Y$ such that  $\phi \sigma_1 = \sigma_2 \phi$.  

Fix a closed oriented surface $F$ of genus $g$ and let $\tau$ be an orientation reversing 
involution of $F$.  The fixed point set $\Fix(\tau)$ of $\tau$ is a collection (possibly empty)
of $m$ disjoint simple closed curves which separates $F$ into $\epsilon$ $=$ $1$ or $2$ 
connected components.  If $\tau_1$ is another such involution, $\tau$ is conjugate to $\tau_1$ 
if and only if $m$ $=$ $m_1$ and $\epsilon$ $=$ $\epsilon_1$ (compare~\cite{S}).  
There are therefore $[(3g+4)/2]$ conjugacy classes of orientation reversing 
involutions on $F$.  $\Fix(\tau)$ is non-separating in $g+1$ cases, while in the remaining 
cases $\Fix(\tau)$ separates $F$ into 
two components, each of genus $(g-m+1)/2$.

Choose $2s$ distinct, ordered points $\{q_i\}$ on $F$ and let 
$\tau$ be an orientation reversing involution of $F$ which 
acts without fixed
points on the $\{q_i\}$.   We define the Teichm\"{u}ller space 
$\cT_g^s (\tau)$ to be the space of all
isomorphism classes of quadruples $(X,\sigma,\{p_i\},f)$ where 
$(X,\sigma)$ is a real algebraic 
curve, $\{p_i\}$ is an ordered collection of $2s$ distinct 
points of $X$ on which $\sigma$ acts freely and 
$$
f: (X,\{p_i\}) \rightarrow (F,\{q_i\})
$$ 
is a homeomorphism preserving the ordering of the points such that $\tau f$ $=$ $f \sigma$.
($F$ is called a {\em marking} of $(X,\sigma,\{p_i\})$.)  
The quadruples $(X,\sigma,\{p_i\},f)$  and $(Y,\mu,\{r_i\},g)$ are isomorphic if 
there exists an isomorphism 
$\phi$ from $(X,\sigma)$ to $(Y,\mu)$ with $\phi(p_i)$ $=$ $r_i$ for each $i$ such 
that $g \phi$ is homotopic to $f \rel \{p_i\}$.   

Since the automorphism group of any Riemann surface $X$ is finite, there
are only finitely many antiholomorphic involutions $\sigma$ of $X$ (any
two differ by an automorphism).  The natural map $\cT_g^s (\tau)$ $\rightarrow$
$\cT_g^{2s}$ which forgets the involution is therefore finite-to-one.  
Topologize $\cT_g^s (\tau)$ via this map.  By choosing an 
appropriate pair of pants decomposition of $F$, it is straightforward 
to use Fenchel-Nielsen coordinates to show that $\cT_g^s (\tau)$ 
is homeomorphic to Euclidean space of dimension $3g-3+2s$. 
(Here $3g-3$ dimensions come 
from the isomorphism class of the real curve
while $2s$ more dimensions come from the locations of 
$s$ of the points and these determine that of the others.) 

It can be shown that when $(X,\sigma,\{p_i\},f)$  is isomorphic 
to $(Y,\mu,\{r_i\},g)$ using the definition above, the homotopy 
from $g \phi$ to $f$ can be replaced by an isotopy $H_t$ such 
that at each time $t$, $H_t(p_i)$ $=$ $q_i$ and $H_t \sigma$ 
$=$ $\tau H_t$.   Therefore, if we let 
$$
\pi_X : X  \rightarrow  X/\sigma
$$
and
$$
\pi_F : F  \rightarrow  F/\tau
$$

\noindent be the quotient maps,  
$\mybar{p}_i$ $=$ $\pi(p_i)$ and $\mybar{q}_i$ $=$ $\pi(q_i)$,
the isomorphism class of $(X,\sigma,\{p_i\},f)$ is determined 
uniquely by $(X/\sigma,\{\mybar{p}_i\})$ and the homotopy 
class  $\rel\{\mybar{p}_i\}$ 
of the homeomorphism of pairs 
$(X/\sigma,\{\mybar{p}_i\}) \rightarrow  (F/\tau,\{\mybar{q}_i\})$ 
induced by $f$.  Thus $\cT_g^s (\tau)$ may be thought of as the
Teichm\"{u}ller space of non-orientable surfaces of genus $g$ (Euler 
characteristic $1-g$) with $s$ ordered points and
a local orientation of the surface at each point.

Next, recall that the mapping class group $\Gamma_g^{2s}$ is 
the group of isotopy classes $\rel \{q_i\}$ of self-homeomorphisms 
of $F$ that fix each $q_i$.  We define the mapping class 
group $\Gamma_g^s(\tau)$ to be the subgroup of 
$\Gamma_g^{2s}$ consisting of all mapping classes which 
admit a representative 
which commutes with $\tau$.  It can be shown that two 
such representatives are homotopic $\rel\{q_i\}$ if and 
only if they are isotopic via an isotopy which commutes 
with $\tau$ at each time.  Since every self-homeomorphism 
of $F/\tau$ has a unique orientation preserving lift to $F$ 
which commutes with $\tau$ 
(the two lifts differ by $\tau$ which is orientation 
reversing), $\Gamma_g^s(\tau)$ may be identified 
with the group of isotopy classes ($\rel\{\mybar{q}_i\}$) 
of self-homeomorphisms of $F/\tau$ 
which preserve each point $\mybar{q}_i$ {\em and} preserve 
a local orientation of $F/\tau$ at each $\pi(q_i)$.  
(Notice that unlike the usual mapping class group of a 
surface with boundary, we 
allow boundary components of $F/\tau$ to be permuted 
and we 
do {\em not} require that homeomorphisms or isotopies of 
homeomorphisms fix boundary curves pointwise.)  
The group $\Gamma_g^s(\tau)$ acts properly discontinuously on 
$\cT_g^s (\tau)$ via the usual action of $\Gamma_g^{2s}$ on
$\cT_g^{2s}$.  The quotient will be denoted $\cM_g^s(\tau)$ and 
called the moduli space of $s$-pointed real algebraic curves.  
It is an orbifold (V-manifold) of real dimension $ 3g-3+2s$  
and, like the moduli space of complex curves, 
it has a finite cover which is an ordinary manifold.   
When two involutions of $F$ are conjugate, there is a 
natural identification of one moduli space with the other, 
so there are $[(3g+4)/2]$ such moduli spaces for each $g$ 
and $s$.

\subsection{The orbifold Euler characteristics}
The Euler characteristic of an orbifold $\cM$ which 
has a finite manifold branched covering $\widetilde{\cM}$ of degree $d$ is 
defined to be

$$
\chi(\cM) = \frac{\chi(\widetilde{\cM})}{d}
$$
where $\chi(\widetilde{\cM})$ is the ordinary Euler characteristic 
of $\widetilde{\cM}$.  This number is independent of the choice 
of $\widetilde{\cM}$ because Euler characteristics multiply by degree 
for unbranched coverings.  Recall that the Euler 
characteristic of a group $G$ which is
virtually torsion free is defined similarly, 

$$\chi(G) = \frac{\chi(\widetilde{G})}{d} $$
where $\widetilde{G}$ is a torsion free subgroup of finite index 
$d$ in $G$.
Since $\cT_g^s (\tau)$  is contractible and $\Gamma_g^s(\tau)$ 
acts properly discontinuously on $\cT_g^s (\tau)$ with 
quotient $\cM_g^s (\tau)$, 

$$
\chi(\cM_g^s (\tau)) = \chi(\Gamma_g^s(\tau)).
$$


\subsection{The Euler characteristic for the real case}

Let $\tau_0$ be an orientation reversing involution of $F$ without 
fixed points. 
Let $B_g$ be the $g$th Bernoulli number, defined by the exponential generating series
\begin{equation}\label{berdef}
B(t)=\frac{t}{e^t-1}=\sum_{j\ge0}B_j\frac{t^j}{j!}.
\end{equation}
Then the  virtual Euler characteristic for the real case is given by the
following result.

\begin{theorem}\label{main}
 For $g \geq 1, s \geq 0 , g+s > 1$, 
$$ \chi(\cM_g^s(\tau_0)) 
= (-2)^{s-1} (1-2^{g-1}) \frac{(g+s-2)!}{g!} B_g.$$
\end{theorem}

>From this theorem it is clear that 
$\chi(\cM_g^s(\tau_0)) = 0$
when $g$ is odd and $g+s > 1$ since $B_{2n+1}$ $=$ $0$ for $n \ge 1$.
The group $\Gamma_1(\tau_0)$ is the subgroup of 
$\Gamma_1$ consisting of all maps which commute 
with $\tau_0$.  Now $\Gamma_1$ is isomorphic to 
$SL_2(\integ)$ since every element of $\Gamma_1$ is determined
by the map it induces on $\pi_1(T^2)$ = $\integ \oplus\integ$.
The map $\tau_0$ induces 
$$ \left( \begin{array}{rr}
            -1 & 0 \\
             0 & 1
          \end{array}  \right)
$$
on $\pi_1(T^2)$, so $\Gamma_1(\tau_0)$ is simply  $\{\pm I\}$.
Since $\Gamma_0$ and $\Gamma_0^2$ are trivial, so are
$\Gamma_0(\tau_0)$ and $\Gamma_0^1(\tau_0)$.  
Also, $\pi_0^2$ is infinite cyclic, isomorphic to
$\Gamma_0^2(\tau_0)$. 
Therefore the above theorem may be completed with the statements:
\begin{eqnarray*}
\chi(\cM_1^0(\tau_0)) = \tsh,\quad
\chi(\cM_0^s(\tau_0)) = 1 \quad\mbox{for $s=0$ or $1$},\quad
\chi(\cM_0^s(\tau_0)) = 0 \quad\mbox{for all $s\ge2$}.
\end{eqnarray*}

For all $g$ and $s$ with $g+s \ge 1$ there are exact sequences 
$$1 \rightarrow 
\pi_g^s \rightarrow  
\Gamma_g^{s+1}(\tau_0)
\rightarrow 
\Gamma_g^s(\tau_0)
\rightarrow 1$$
\noindent where $\pi_g^s$ is the fundamental group of the orientable 
surface of genus $g$ with $2s$ points removed.  The map 

$$\Gamma_g^{s+1}(\tau_0) \rightarrow \Gamma_g^s(\tau_0)$$ 
\noindent is obtained by forgetting one of the pairs of points, say 
$(q_{2s+1},q_{2s+2})$, which are interchanged by
$\tau_0$.  Its kernel is $\pi_g^s$ because the isotopy to
the identity of an element of $\Gamma_g^{s+1}(\tau_0)$ 
which fixes $\{q_1,...,q_{2s}\}$ creates and is
created by moving the point $q_{2s+1}$ along a loop 
representing an element of the fundamental group of 
$F - \{q_1,...,q_{2s}\}$.
($F$ occurs here instead of $F/\tau_0$ because the points ${q_i}$ 
are not permuted.  Equivalently, the self-maps of $F/\tau$ 
which we consider are all locally orientable at the points 
$\pi(q_i)$, so only orientation preserving loops lie in 
$\Gamma_g^{s+1}(\tau_0)$.)   

Now, given any short exact sequence of groups
$$1 \rightarrow A \rightarrow B \rightarrow C \rightarrow 1$$ 
we have
$$\chi(B) = \chi(A) \chi(C).$$

This means that for each $g$, Theorem~\ref{main} only needs to be proven for one
value of $s$.  We will, in fact, show it for all $s \ge 1$, as this is the
more natural approach.  The short exact sequence is necessary, however, to
establish the case where $s = 0$.  

For the involutions $\tau$ which have fixed points, let us consider separately the cases 
where the fixed point sets do or do not separate $F$.  If $\Fix(\tau)$ consists of $m$ disjoint 
simple closed curves which together do not separate $F$, then $F/\tau$ is a non-orientable 
surface with $n$ boundary components.  The Dehn 
twist on a boundary curve of $F/\tau$ is isotopic to the identity 
since it lifts to two Dehn twists on isotopic curves with
opposite twist directions.  
In addition, elements can permute the components of 
$\Fix(\tau)$.  Therefore $\Gamma_{g-m}^{s+m} (\tau_0)$ may 
be identified with a subgroup of index $m!$ in the mapping 
class group of $(F/\tau,\{\pi(q_i)\})$.  
The Euler characteristic of the former is given by
Theorem~\ref{main} above, that of the latter is obtained by
dividing by $m!$

If $\Fix(\tau)$ consists of $m$ disjoint simple closed 
curves which together do separate $F$, then 
$F/\tau$ is an orientable surface with $m$ 
boundary components of genus $h$, where 
$g = 2h + m - 1$.  In this case the ordinary 
mapping class group $\Gamma_h^{s+m}$ may be 
identified with a subgroup of index $m!$ in the 
mapping class group of $(F/\tau,\{\pi(q_i)\})$.  
The Euler characteristic of $\Gamma_h^{s+m}$ 
was computed in~\cite{HZ}.  Combining all 
of these results we have the following general statement.

\begin{corollary}
If $\Fix(\tau)$ is non-separating and consists of $m$
 simple closed curves then
\begin{eqnarray*}
\chi(\Gamma_g^s (\tau)) &=&
(-2)^{s+m-1} (1-2^{g-m-1}) \frac{(g+s-2)!}{m! (g-m)!} B_{g-m}.
\end{eqnarray*}

If $\Fix(\tau)$ is separating and consists of $m$ simple 
closed curves then $g-m+1$ must be even and
\begin{eqnarray*}
\chi(\Gamma_g^s (\tau)) =
(-1)^{s+m} \frac{(g-m+s-2)!}{m! (g-m+1) (g-m-1)!} B_{g-m+1}
\end{eqnarray*}
\end{corollary}


\subsection{Moduli spaces and polygonal identifications}

We begin the proof of Theorem~\ref{main} by first expressing $\chi(\cM_g^s(\tau_0))$ in
terms of pairings of sides of polygons.
Recall~\cite{Str} that for 
every $k$-pointed 
Riemann surface $(X,\{p_1,\ldots,p_k\})$ and every
collection of weights $(t_1,\ldots,t_k)$ with $t_i > 0$ and $t_1+\ldots+t_k$ 
$=$ $1$, there 
is a unique quadratic differential which has double poles at the points ${p_i}$ of 
type ${-dz^2}/{z^2}$, no other poles, closed real trajectories 
and connected singular trajectories 
with the distance (in the norm defined by the quadratic 
differential) from 
$p_i$ to the singular trajectory equal to $t_i$.
This concept was first used in~\cite{Harer2} (see also~\cite{Harer1})
to describe cells for the moduli space of curves.  
The singular trajectories of the differential then decompose
$X$ into $k$ disks which, after scaling $X$,
can be taken to be radius $t_i$, centered at each $p_i$.  This
exhibits $(X,\{p_1,\ldots,p_k\})$ as obtained in a unique manner by
identifying the sides of $k$ disks $D_i$ in
$\comp$
(compare~\cite{Harer1} and~\cite{Harer2}).  Fixing an identification pattern
with a marking and varying the edge lengths determines
an open cell in $\cT_g^k \times \Delta^{k-1}$, where $\Delta^{k-1}$ is a
$k-1$ simplex.

Suppose now that $X$ admits a fixed point free involution $\tau_0$
which interchanges the points $p_{2i-1}$ and $p_{2i}$ for $i=1,\ldots, s$.
It follows directly from Strebel's work that the real and
imaginary trajectories of the differential are invariant
under $\tau_0$.  (This fact was first 
pointed out to us by Sepp\"{a}l\"{a}.)
This means that $\cT_g^s(\tau_0) \times \Delta^{s-1}$ 
admits an ideal triangulation with cells determined
as follows: Choose positive integers $n_1,\ldots,n_s$ with
$n_1+\ldots +n_s$ $=$ $2n$. 
For each $i$ take two copies $P_i^+$ and $P_i^-$ of a polygon 
with $n_i$ sides.  Also take a pairing $\omega$ of the sides of 
$ \cup_{1\leq i \leq s} P_i^+ \cup P_i^-$ 
which gives an oriented, connected surface of genus $g$ and is compatible
with the interchange maps $P_i^+ \longleftrightarrow P_i^-$.  These maps
together induce an involution $\sigma_\omega$ of 
$ (\cup_{1\leq i \leq s} \partial P_i^+ \cup \partial P_i^-)/\omega$
We assume that $\sigma_\omega$ is orientation reversing.  Finally, this 
 pairing 
(as in~\cite{Harer2} and~\cite{HZ}) must satisfy the requirement 
that the valence of each vertex of the boundary graph 
$ (\cup_{1\leq i \leq s} \partial P_i^+ \cup \partial P_i^-)/\omega$ 
is at least $3$.  Finally, choose a marking $f_{\omega}$
of the result that identifies the center of $P_i^+$ with $q_{2i-1}$
and the center of $P_i^-$ with $q_{2i}$.  Varying the edge lengths
determines an open cell whose dimension is $n-1$.
Orbits of cells are 
determined by the identification alone (without the marking).  
This gives an {\em ideal} triangulation in the following sense:
Each marked identification determines an open cell, corresponding 
to positive lengths on each edge.  By allowing an edge length to
go to $0$, we obtain a cell of one lower dimension as long as the
resulting identification still gives a surface of genus $g$.  If,
however, the result has smaller genus, the face lies outside of 
$\cT_g^s(\tau_0) \times \Delta^{s-1}$.  

By forming the quotient of $X_{\omega}$ by $\tau_0$, the identification
$\omega$ induces an identification $\omega^{\prime}$ of the union of
$s$ polygons $P_i$ having $n_i$ edges to obtain a connected,
non-orientable surface $X_{\omega}^{\prime}$ of Euler characteristic 
$1-g$.  
(The boundary graph has $n$ edges and $1+n-g-s$ vertices and
the same valence condition as above, since $\tau_0$ is fixed point
free.)  The marking $f_{\omega}$ induces a marking $f_{\omega}^{\prime}$
of $X_{\omega}^{\prime}$ which pairs the ordered center points 
$\{\mybar{p}_1,\ldots,\mybar{p}_s\}$ with $\{\mybar{q}_1,\ldots,\mybar{q_s}\}$.

Conversely, given an identification ${\omega}^{\prime}$ of the edges of 
$ \cup_{1\leq i \leq s}  P_i $ to obtain a non-orientable surface 
$X_{\omega}^{\prime}$ of Euler characteristic $1-g$, with special points
$\{\mybar{p}_1,\ldots,\mybar{p}_s\}$, and given a marking $f_{\omega}^{\prime}$
of $X_{\omega}^{\prime}$ with $f_{\omega}^{\prime}(\mybar{p}_i)$ $=$
$\mybar{q}_i$, there are $2^{s-1}$ corresponding choices of the pairing 
$\omega$ and marking $f_\omega$.  To see this, visualize the disks $P_i^+$
and $P_i^-$ in a common plane so that assigning $U$ for an up normal
or $D$ for a down normal orients them.  Fix $U$ as the orientation for 
$P_1^+$ and $D$ for $P_1^-$.  Arbitrarily assign $U$ or $D$ to $P_i^+$
and the opposite to $P_i^-$ for $i = 2, \ldots, s$.  Now there is one and 
only one identification $\omega$ lifting $\omega^\prime$ which gives a 
surface on which the chosen orientations agree:

Think of the disks $P_i$ as also lying in a common plane with a $U$ 
orientation on each.  Let $e_1$ be an edge of $P_i$ and $e_2$ an edge of
$P_j$ (perhaps $i=j$) and suppose that $\omega^\prime$ pairs $e_1$ with
$e_2$.  Let $e_1^U$ be the edge corresponding to $E_1$ in either
$P_i^+$ or $P_i^-$, whichever is labeled $U$.  Similarly define 
$e_1^D$, $e_2^U$ and $e_2^D$.  If $\omega^\prime$ pairs $e_1$ and $e_2$ 
so that the orientations on $P_i$ and $P_j$ agree, then $e_1^U$ is
paired with $e_2^U$ and $e_1^D$ with $e_2^D$.  Otherwise they switch.
(Notice that we could have performed the entire process by first starting
with any $P_i^+$ as either up or down and we would have the same $2^{s-1}$
possibilities.)  The marking $f_{\omega}^{\prime}$ lifts to a unique 
$f_{\omega}$ once we specify which point of $\pi^{-1}(\mybar{q}_1)$ to call
$q_1$.  Therefore, given $({\omega}^{\prime},f_{\omega}^{\prime})$, there 
are exactly $2^{s-1}$ distinct lifts $({\omega},f_{\omega})$.  Since 
elements of the 
mapping class group $\Gamma_g^s(\tau_0)$ do not permute the points $\{q_i\}$,
there are similarly $2^{s-1}$ distinct $\omega$.

We now consider polygonal identifications induced
by $\omega^{\prime}$.
Define $\lambda_g^N(n_1,\ldots,n_s)$ to be the number of identifications 
$\omega$ of the edges of 
$ \cup_{1\leq i \leq s} P_i$ 
that give a non-orientable surface of genus $g$ (Euler characteristic
$1-g$) and so that the  valence of each vertex of the 
boundary graph 
$ (\cup_{1\leq i \leq s} \partial P_i )/\omega$
is at least $3$.  Here we are assuming that 
an initial edge of each $P_i$ is fixed so that, for example 
when $s$ $=$ $1$, 
the identifications $a a b b $ and $a b b a$ would be counted 
as different ways of constructing a Klein bottle.
As an illustration,
\begin{equation}\label{klein}
\lambda_1^N(2)=4, 
\end{equation}
corresponding to the four patterns:
$$a\,a\,b\,b,\quad a\,b\,b\,a,\quad a\,b\,a\, b^{-1},\quad a\,b\,a^{-1}\,b.$$
Notice that if the surface has Euler characteristic $1-g$ then the number of vertices 
of the boundary graph is $1+n-g-s$.  

Similarly, we define 
$\lambda_g^O(n_1,\ldots,n_s)$ to be the number of identifications 
of the edges of 
$ \cup_{1\leq i \leq s} P_i$ 
that give a {\em orientable} surface of genus $g$ (Euler characteristic
$2-2g$) with the same valence condition.  Finally, define 
$\lambda_g(n_1,\ldots,n_s)$ to be the number of identifications 
of the edges of 
$ \cup_{1\leq i \leq s} P_i$ 
that give a surface of Euler characteristic
$1-g$, {\em orientable or not}, with the same valence condition. 

Now set
$$
\lambda_g^s(n) = \sum_{n_1+\cdots+n_s=2n}  \lambda_g(n_1,\ldots,n_s),
$$
and let $\lambda_g^{s,N}(n)$ and $\lambda_{(g+1)/2}^{s,O}(n)$ be the sums, over the same
range, with summands $\lambda_g^N(n_1,\ldots,n_s)$ and $\lambda_{(g+1)/2}^O(n_1,\ldots,n_s),$
respectively.
Finally, let
\begin{eqnarray}\label{lamchi}
\Lambda_g^s = \sum_{n=g+s}^{3g+3s-3} \frac{(-1)^{n-s}}{2n} \lambda_g^s(n)
\end{eqnarray}
and let $\Lambda_g^{s,N}$ and $\Lambda_{(g+1)/2}^{s,O}$ be the analogous sums
in which $\lambda_g^s(n)$ is replaced by $\lambda_g^{s,N}(n)$ and $\lambda_{(g+1)/2}^{s,O}(n),$
respectively.

Next we express $\chi(\cM_g^s(\tau_0))$ in
terms of one of the combinatorial sums defined above.

\begin{lemma}\label{T23}
$$\chi(\cM_g^s(\tau_0)) = 2^{s-1} \Lambda_g^{s,N}.$$
\end{lemma}
\proof   Recall that $\cM_g^s(\tau_0))$ is an orbifold and by
definition, each point of an orbifold $\cM$ has a 
neighborhood that is locally modeled on $\reals^{\dim(\cM)}$ modulo a
finite group.  An orbifold triangulation of $\cM$ is a triangulation
by simplices $\omega$
which have the property that the finite group associated to 
each point of the interior of $\omega$ is the same.  
If we call this group
$G(\omega)$, the Euler characteristic on $\cM$ is computed with the formula:

$$\chi(\cM) = \sum_{\omega} \frac{(-1)^{\dim(\omega)}}{o(G(\omega))}.$$

To compute the Euler characteristic of $\cM_g^s(\tau_0)$, consider the
ideal triangulation of $\cM_g^s(\tau_0) \times \Delta^{s-1}$ described by
identification patterns.  We cannot apply the formula above directly because 
cells of the triangulation must be closed.  Instead we use the dual complex
described in~\cite{Harer2}: $\cM_g^s(\tau_0) \times \Delta^{s-1}$ retracts onto
a spine complex $Y$ of dimension $2g-3+3s$ which has a $k$-cell for each 
identification pattern $\omega$ with $n-1 = 3g-3s-4-k$.  
The cells are the identification patterns and the 
local groups are the symmetries of the
configuration which determines the cell.  Since these must
fix the points $q_i$, they are all cyclic. 
Counting each
identification pattern of (with no initial edges) weighted 
by the reciprocal of the order
of the cyclic symmetry group is the same as counting each
pattern with a choice of initial edge weighted by 
${1}/{2n}$
(compare~\cite{HZ}).
The result follows. \hfill\qed

\medskip
Now clearly 
\begin{eqnarray}\label{xidiff} 
\Lambda_g^s = \Lambda_g^{s,N} + \Lambda_{(g+1)/{2}}^{s,O},
\end{eqnarray}
where the second term on the right is naturally zero when $g$ is even.
Then Lemma~\ref{T23} and~(\ref{xidiff}) enable us to establish of Theorem~\ref{main} result
by determining $\Lambda_g^s$ and $\Lambda_{(g+1)/{2}}^{s,O}.$ 
This will be carried out in the next section, where
we actually determine a parametrized quantity that specializes
to each of $\Lambda_g^s$ and $\Lambda_{(g+1)/{2}}^{s,O}$ at different values of the
parameter.

Previously Harer and Zagier~\cite{HZ}  showed that
$\chi(\cM_g^s)= \Lambda_{(g+1)/{2}}^{s,O}$ and thus 
determined that
\begin{eqnarray*}
\chi(\cM_g^s) = (-1)^s\frac{(g+s-2)!}{(g+1)(g-1)!}B_{g+1},
\end{eqnarray*} 
when $g$ is odd
(of course, $\chi(\cM_g^s)=0$ for $g$ even).

\section{Matrix models and the counting series for embedded graphs}\label{Sge}

\subsection{Polygonal identifications and graph embeddings}

To determine $\Lambda_g^s$ and $\Lambda_{(g+1)/2}^{s,O}$ we first
identify $\lambda_g^s(n)$ and $\lambda^{s,O}_{(g+1)/2}(n)$ as solutions
of equivalent enumerative questions for combinatorial maps,
since appropriate forms for their generating series are
already available, from current work in algebraic combinatorics.
A {\em map} is an embedded graph in a surface with the property
that each complementary region is a 2-cell.  In topology one says
the graph fills the surface and in combinatorics ones calls this
a 2-cell embedding.
The graph divides
the surface into regions called {\em faces}. 
In an orientable surface a map is {\em rooted} by distinguishing an edge, 
and an end of this edge. In nonorientable surfaces
an edge is distinguished together 
with an end of this edge, and a side of this end of the edge. Rooted maps have 
only the trivial automorphism. 
From this point of view, an identification $\omega$ of the edges
of $P_n$ with boundary graph $(\partial P_n)/\omega$ corresponds
precisely to a two-cell embedding of the boundary graph with
a single face; the rooting allows us to uniquely recover
the labels of the sides and the orientation of $P_n$. Thus
we have the following combinatorial identifications:

\begin{itemize}
\item $\lambda^{s,O}_g(n)$ is $s!$ times
the number of rooted maps with s faces, $n$ edges and no vertices
of valences $1$ or $2,$ in orientable  surfaces of genus $g$ (Euler
characteristic $2-2g$),
\item $\lambda_g^s(n)$ is $s!$ times
the number of rooted maps with s faces, $n$ edges and no vertices
of valences $1$ or $2,$ in locally orientable surfaces of Euler
characteristic $1-g$.
\end{itemize}

For $\bfi=(i_1,i_2,\ldots),$  let $m(\bfi,j,n)$ and $m^O(\bfi,j,n)$
be, respectively, the numbers of rooted maps in locally orientable and orientable
surfaces, with $n$ edges, $j$ faces and $i_k$ vertices of valence $k$ for $k\ge1$.
Thus $\sum_{k\geq 1}k i_k
=2n$, and $\sum_{k\geq 1}i_k =\vert V\vert$, the number of vertices.
But from the Euler-Poincar\'{e} theorem, for the maps counted
by $\lambda_{(g+1)/2}^O(n)$ and $\lambda_g(n)$ we
have $\vert F\vert -\vert E\vert
+\vert V\vert =s-n+\vert V\vert =1-g$ so $\vert V\vert =n-g-s+1$.
Thus we can write
\begin{equation}\label{reln}
\lambda_g^s(n) = s! \sum_\bfi m(\bfi ,s,n),\quad\quad
\lambda_{(g+1)/2}^{s,O}(n) = s! \sum_\bfi m^O(\bfi ,s,n),
\end{equation}
where, in both cases the sum is over $\bfi$ such that
\begin{equation}\label{cond}
i_1=i_2=0,\quad\quad\sum_{k\geq 3}k i_k =2n,\quad\quad
\sum_{k\geq 3}i_k=n-g-s+1.
\end{equation}

The vector $\bfi$ is called the {\em vertex distribution} of such maps.
Let the  generating series for maps with respect to vertex distribution be
\begin{equation}\label{ee2a}
M(\bfy,x,z) = \sum_{\bfi,j,n}m(\bfi,j,n)\bfy^\bfi x^jz^n,\;\;\;
M^O(\bfy,x,z) = \sum_{\bfi,j,n}m^O(\bfi,j,n)\bfy^\bfi x^jz^n,
\end{equation}
where $\bfy = (y_1,y_2,\ldots)$, and $\bfy^\bfi =\prod_{k\geq 1} y_k^{i_k}$.

\subsection{A pair of matrix integrals}

These two series have the following integral representations,
as given in Jackson~\cite{JI}, for $M^O$, and Goulden
and Jackson~\cite{GJI}, for $M$:
\begin{eqnarray*}
M^O(\bfy,N,z) &=&2z\frac{\partial}{\partial z}\log\left( \frac
{\displaystyle{\int_{\cV_N}e^{\sum_{k\ge1}\frac{1}{k}y_k\sqrt{z}^k\tr\bfM^k} e^{-\tsh\tr\bfM^2}d\bfM}
}{\displaystyle{\int_{\cV_N} e^{-\tsh\tr\bfM^2}d\bfM}}\right),\\
M(\bfy,N,z) &=&4z\frac{\partial}{\partial z}\log\left( \frac
{\displaystyle{\int_{\cW_N}e^{\tsh\sum_{k\ge1}\frac{1}{k}y_k\sqrt{z}^k\tr\bfM^k}
e^{-\tsqu\tr\bfM^2}d\bfM}
}{\displaystyle{\int_{\cW_N} e^{-\tsqu\tr\bfM^2}d\bfM}}\right),
\end{eqnarray*}
where $\cV_N$ and $\cW_N$ are, respectively, the vector spaces of  Hermitian complex and
real symmetric  $N\times N$  matrices.
The reader who wishes to have greater detail is directed to the above
sources, and to the discussion in Section~\ref{Scrvec} of this paper.

The integrals over $\cV_N$ and $\cW_N$ can be transformed by the Weyl integration theorems~\cite{HE} to
the following
integrals over the reals.  
For a positive integer $N,$ let $\lambda=(\lambda_1,\ldots,\lambda_N),$ let $V(\lambda)$ be the
Vandermonde determinant $\prod_{1\le i<j\le N}(\lambda_j-\lambda_i)$, 
 let $p_k(\lambda)=\lambda_1^k+\cdots+\lambda_N^k,$ the
power sum symmetric function of degree $k,$ for $k\geq 1,$
and $d\lambda=d\lambda_1\cdots d\lambda_N$ :
\begin{eqnarray}\label{orint}
M^O(\bfy,N,z)=2z\frac{\partial}{\partial z}\log \left(
\frac{
\displaystyle{\int_{\reals^N} \vert V(\lambda)\vert^{2}
e^{\sum_{k\ge1}\frac{1}{k}\sqrt{z}^ky_kp_k(\lambda) }
e^{-\frac{1}{2}\, p_2(\lambda)}d\lambda
}
}{
\displaystyle{\int_{\reals^N} \vert V(\lambda)\vert^{2} 
e^{-\frac{1}{2}\, p_2(\lambda)}d\lambda 
}
}\right),
\end{eqnarray}
\begin{eqnarray}\label{locint}
M(\bfy,N,z)=4z\frac{\partial}{\partial z}\log \left(
\frac{
\displaystyle{\int_{\reals^N} \vert V(\lambda)\vert
e^{\tsh\sum_{k\ge1}\frac{1}{k}\sqrt{z}^ky_kp_k(\lambda) }
e^{-\frac{1}{4}\, p_2(\lambda)}d\lambda
}
}{
\displaystyle{\int_{\reals^N} \vert V(\lambda)\vert 
e^{-\frac{1}{4}\, p_2(\lambda)}d\lambda 
}
}\right).
\end{eqnarray}

A joint generalization of these is afforded by
\begin{eqnarray}\label{ee4}
M_\gamma(\bfy,x,z)=\frac{2}{\gamma}z\frac{\partial}{\partial z}
\log Z_\gamma(\bfy,x,z),
\end{eqnarray}
where 
\begin{eqnarray}\label{ee5}
Z_\gamma (\bfy,N,z) =
\frac{
\displaystyle{\int_{\reals^N} \vert V(\lambda)\vert^{2\gamma}
e^{\gamma\sum_{k\ge1}\frac{1}{k}\sqrt{z}^ky_kp_k(\lambda) }
e^{-\frac{1}{2}\gamma\, p_2(\lambda)}d\lambda
}
}{
\displaystyle{\int_{\reals^N} \vert V(\lambda)\vert^{2\gamma} 
e^{-\frac{1}{2}\gamma\, p_2(\lambda)}d\lambda 
}
},
\end{eqnarray}
since
\begin{equation}\label{specgam}
M(\bfy,N,z)=M_{\tsh}(\bfy,N,z),\;\;\; M^O(\bfy,N,z)=M_1(\bfy,N,z).
\end{equation}

Thus both $M$ and $M^O$ can be obtained by specializing the
parameter $\gamma$ in $M_\gamma$, and
it is the series $M_\gamma$ that we therefore study in detail.
We now make a few comments about the interpretation of the above
series. First, the numerator of~(\ref{ee5}) is an even function
of $\sqrt{z}$ (to see this, consider the
simultaneous substitutions $\lambda_j \rightarrow -\lambda_j ,\, j=1,\ldots ,N$),
so the occurrence of $\sqrt{z}$ is of no
concern; $M_\gamma$ really is a power series
in $z$. Second, the second argument in these map series is
an indeterminate, $x$, yet in each case the integral representation
evaluates the series at any positive integer value $N$ (which is the
dimension of the region of integration). However,
the number of rooted maps with $n$ edges is finite for each $n$, so
 $M_\gamma$ can be viewed as a power
series in $z$ with coefficients that are
polynomials in $y_1,\ldots ,y_{2n},x$. Thus, since this holds
for infinitely many $N$, a polynomiality argument
allows us to obtain $M_\gamma(\bfy,x,z)$ 
from $M_\gamma(\bfy,N,z)$, in principle, by deriving
an explicit presentation of $M_\gamma(\bfy,N,z)$ as a power series in $z$,
and formally replacing $N$ by $x$ in the resulting (polynomial) coefficients
of each power of $z$. 

In the following result, we obtain preliminary
expressions for $\Lambda^s_g$ and $\Lambda_{(g+1)/2}^{s,O}$, as
specializations of the parameter in a quantity
$\xi^s_g(\gamma )$, which is therefore referred to as the
{\em parametrized Euler characteristic}.
Here we use  the notation $[A]B$ for the
coefficient of $A$ in the expansion of $B$.
The result gives $\xi^s_g(\gamma)$ as a particular coefficient of $M_\gamma$ evaluated
at transformed arguments.

\begin{proposition}\label{p1} Let
\begin{equation}\label{defH}
{\grgen}^s_g (\gamma)= s!(-1)^s[x^s t^{g+s-1}]\,
\frac{1}{\gamma}\log W_\gamma(x,t),
\end{equation}
where
\begin{eqnarray}\label{ee6}
W_\gamma(N,t)=
\frac{
\displaystyle{
\int_{\reals^N} \vert V(\lambda)\vert^{2\gamma}
\left(\prod_{j=1}^N e^{-i\gamma\lambda_j/ \sqrt{t}}
(1-i\sqrt{t}\lambda_j)^{-\gamma/t}\right)d\lambda
}
}{
\displaystyle{
\int_{\reals^N} \vert V(\lambda)\vert^{2\gamma}e^{-\gamma p_2(\lambda)/2}d\lambda
}
}.
\end{eqnarray}
Then, for $s,g\geq1,$
\begin{equation}\label{specpar}
{\Lambda}^s_g = {\grgen}^s_g (\tsh ),\quad\quad{\Lambda}_{(g+1)/2}^{s,O} =
{\grgen}^s_g (1).
\end{equation}
\end{proposition}
\proof
Let ${\bf u}(t) =(u_1,u_2,\ldots)$, where $u_1=u_2=0,
\; u_k=-(i\sqrt{t})^{k-2},
k\ge3$ (here $i$ is $\sqrt{-1}$, and is not
an index), and let $\Psi$ be the operator whose action is defined by
$$\Psi f(\bfy,x,z) ={\tsh}\int_0^1 f({\bf u}(t),x,z)\frac{dz}{ z}.$$

Now, for $\bfi$ satisfying the conditions in~(\ref{cond}), under
the specialization $\bfy = {\bf u}(t)$, we obtain
\begin{eqnarray*}
\bfy^\bfi =
(-1)^{n-g-s+1+\tsh\{2n-2(n-g-s+1)\}}\,
t^{\tsh\{2n-2(n-g-s+1)\}} 
= (-1)^n t^{g+s-1}.
\end{eqnarray*}
Then we obtain
$\Lambda^s_g=s!(-1)^s[x^s t^{g+s-1}]\Psi M$ and
$\Lambda_{(g+1)/2}^{s,O}=s!(-1)^s[x^s t^{g+s-1}]\Psi M^O,$
for $g\geq 1,s\geq 1,$ from~(\ref{lamchi}),~(\ref{reln})
and~(\ref{ee2a}).
But, from~(\ref{ee4}) and~(\ref{ee5}),
$$\Psi M_\gamma = \frac{1}{\gamma}\int_0^1\left( \frac{\partial}{\partial z}
\log Z_\gamma (\bfu (t) ,N,z)\right) dz
= \frac{1}{\gamma}\log  Z_\gamma (\bfu (t),N,1),$$
since $Z_\gamma (\bfu (t),N,0)=1$.
Now let $W_\gamma (N,t)=Z_\gamma (\bfu (t),N,1)$. The result follows
from~(\ref{specgam}) and the simplification
\begin{eqnarray*}
-\tsh p_2(\lambda )+\sum_{k\geq 1} \frac{1}{k}u_k p_k(\lambda )&=&
-\sum_{k\geq 2} \frac{{(i\sqrt{t})^{k-2}} }{ k}\, p_k (\lambda )\\
&=&\frac{1}{i\sqrt{t}}\, p_1(\lambda )+\frac{1}{t}
\sum_{k\geq 1}\frac{{(i\sqrt{t})^k}}{ k}
\, p_k (\lambda )\\
&=&\sum_{j=1}^N \left( \frac{-i\lambda_j }{\sqrt{t}}+\frac{1}{t}\log(1-i\sqrt{t}
\lambda_j )^{-1}\right) .
\end{eqnarray*}
Note that the application of $\Psi$ to $M_\gamma$
is justified,
since the number of maps with no vertices of valences $1$ or $2$ (this
is forced by the specialization $y_1 =y_2 =0$), and $\vert E\vert
-\vert V\vert =g+s-1$ is finite for
each $g$ and $s$ (e.g., these conditions imply
that there are at most $3g+3s-3$ edges, as in the upper bound for
the summation in~(\ref{lamchi})).
Thus $W_\gamma (x,t)$ can
be viewed as a power
series in $t$ with coefficients that are polynomials in $x$.
\hfill\qed

\section{Determination of the parametrized Euler characteristic}\label{Spmr}
The following detailed calculation involving properties of the
gamma function of the parameter $\gamma$ is necessary since we wish
to obtain the parametrized Euler characteristic $\xi^s_g(\gamma)$ as
 an explicit polynomial in $\gamma^{-1}.$ In Section~\ref{Scrvec} we shall show that
$b=\gamma^{-1}-1$ is an indeterminate that is conjectured to be associated with a
combinatorial invariant of rooted maps. 

\subsection{Selberg integration}

We now  determine $W_\gamma (x,t)$ from the
integral representation~(\ref{ee6}) for $W_\gamma (N,t)$
by giving an explicit presentation of $W_\gamma (N,t)$ as a series in $t$,
with coefficients that are polynomials in $N$, and then, appealing
to polynomiality, formally replacing $N$ by $x$.
To evaluate $W_\gamma(N,t)$ we use the following
theorem due to Mehta~\cite{ME}. It is derived in a
similar fashion to Selberg's theorem, but based on an integral of Cauchy
instead of the beta integral
 (note that again $i=\sqrt{-1}$, and is not an index).

\begin{theorem}\label{tsic}
For ${\rm Re}\, a,{\rm Re}\, b,{\rm Re}\,\alpha,{\rm Re}\,\beta>0,
{\rm Re}\, (\alpha+\beta)>1$ and
$$-\frac{1}{N}<{\rm Re}\,\gamma<\min\left( \frac{{\rm Re}\,\alpha}{N-1},
\frac{{\rm Re}\,\beta}{N-1},
\frac{{\rm Re}\, (\alpha+\beta+1)}{2(N-1)}
\right) $$
then
$$
\rintn\vert V(\lambda )\vert^{2\gamma}
\prod_{j=1}^N\left(1+i\frac{\lambda_j}{a}\right)^{-\alpha}
\left(1-i\frac{\lambda_j}{b}\right)^{-\beta}d\lambda
=\kappa
\prod_{j=0}^{N-1}
\frac{\Gamma(\alpha+\beta-(N+j-1)\gamma-1) }{
\Gamma(\alpha-\gamma j)\Gamma(\beta-\gamma j)
}
$$
where $$\kappa=\left(\frac{2\pi\,a^\alpha b^{\,\beta}
}{(a+b)^{(\alpha+\beta)-\gamma(N-1)-1}}\right)^N
\prod_{j=0}^{N-1}\frac{\Gamma(1+(1+j)\gamma) }{\Gamma(1+\gamma)}.$$
\end{theorem}

\medskip
In proving the following result we use Mehta's integration theorem to
evaluate $\log W_\gamma (N,t)$, and
then apply the polynomiality
argument to give an explicit (asymptotic) power series presentation
of $\log W_\gamma (x,t)$.

\begin{theorem}\label{l3}
As an asymptotic expansion in $1/t$ we have
$$\log W_\gamma (x,t)
=-\gamma x\sum_{k\ge1}\frac{ B_{2k}t^{2k-1} }{ 2k(2k-1)}
+\sum_{\delta\ge1}\frac{t^\delta}{ \delta(\delta+1)}
\sum_{r=1}^{\delta+1}\binom{\delta+1}{r} B_{\delta+1-r}$$
$$\times\left(\frac{{x^r (-1)^{\delta+1-r}} }{ \gamma^{\delta-r} }
-\sum_{m=1}^{r+1}\binom{r+1}{m}\frac{B_{r+1-m}}{r+1}
\frac{x^m}{\gamma^{r-m}}(-1)^{\delta-m}
\right).$$
\end{theorem}
\proof From~(\ref{ee6}), $W_\gamma(x,t)=F_\gamma(x,t)/G_\gamma(x,t),$ where
\begin{eqnarray*}
F_\gamma(N,t)=\int_{\reals^N} \vert V(\lambda)\vert^{2\gamma}
\left(\prod_{j=1}^N e^{-i\gamma\lambda_j/ \sqrt{t}}
(1-i\sqrt{t}\lambda_j)^{-\gamma/t}\right)d\lambda,
\end{eqnarray*}
and
\begin{eqnarray*}
G_\gamma(N,t)
=\int_{\reals^N} \vert V(\lambda)\vert^{2\gamma}e^{-\gamma p_2(\lambda)/2}d\lambda.
\end{eqnarray*}
We begin by determining the integrals $F_\gamma(N,t)$ and   $G_\gamma(N,t)$ 
as different specializations of Mehta's integral in Theorem~\ref{tsic} above.

\noindent {To determine $G_\gamma$:}
We can rewrite $G_\gamma$ in the form
\begin{eqnarray*}
G_\gamma(N,t)=\lim_{L\rightarrow\infty}
\rintn\vert V(\lambda )\vert^{2\gamma}
\prod_{j=1}^N\left(1+i\frac{\lambda_j}{\sqrt{L}}\right)^{-\tsh\gamma L}
\left(1-i\frac{\lambda_j}{\sqrt{L}}\right)^{-\tsh\gamma L}d\lambda,
\end{eqnarray*}
so applying Theorem~\ref{tsic}, with $a=b=\sqrt{L}$ and $\alpha=\beta=\gamma L/2$,
gives
\begin{eqnarray*}
G_\gamma=\left(\prod_{j=0}^{N-1}\frac{\Gamma(1+(1+j)\gamma)}{ \Gamma(1+\gamma)}\right)
\lim_{L\rightarrow\infty}
\left( \frac {2\pi \sqrt{L}^{\gamma L}
}{ (2\sqrt{L})^{\gamma L-\gamma(N-1)-1}
}\right)^N
\prod_{j=0}^{N-1}
\frac{\Gamma(\gamma L-(N+j-1)\gamma-1) }{
\Gamma(\tsh\gamma L-\gamma j)^2
}.
\end{eqnarray*}
But the duplication theorem for the gamma function gives
$$
\Gamma(\gamma L-(N+j-1)\gamma-1) =\frac{1}{\sqrt{\pi}}2^{\gamma L-\gamma(N+j-1)-2}
\Gamma(\tsh(\gamma L-(N+j-1)\gamma-1))\Gamma(\tsh(\gamma L-(N+j-1)\gamma)),
$$
and it is straightforward to establish by standard properties of
the gamma function that
\begin{eqnarray*}
\lim_{L\rightarrow\infty}
\frac{\Gamma(\tsh(\gamma L-(N+j-1)\gamma-1))
}{ \Gamma(\tsh\gamma L-\gamma j)\,\,
(\tsh\gamma L)^{-\tsh(N-j-1)\gamma-\tsh}
}=1=
\lim_{L\rightarrow\infty}
\frac{\Gamma(\tsh(\gamma L-(N+j-1)\gamma))
}{\Gamma(\tsh\gamma L-\gamma j)\,\,
(\tsh\gamma L)^{-\tsh(N-j-1)\gamma}
}.
\end{eqnarray*}
It follows that
\begin{eqnarray*}
G_\gamma(N,t)=\left(\frac{\sqrt{2\pi}}{\gamma^{\tsh((N-1)\gamma+1)}}\right)^N
\prod_{j=0}^{N-1}\frac{\Gamma(1+(1+j)\gamma)}{\Gamma(1+\gamma)}.
\end{eqnarray*}

\noindent {To determine $F_\gamma$:}
We can rewrite $F_\gamma$ in the form
\begin{eqnarray*}
F_\gamma(N,t)=\lim_{L\rightarrow\infty}\int_{\reals^N}
\vert V(\lambda)\vert^{2\gamma}
\prod_{j=1}^N
\left(1+i\frac{\sqrt{t}}{ L}\lambda_j\right)^{-\gamma L/t}
(1-i\sqrt{t}\lambda_j)^{-\gamma/t}d\lambda
\end{eqnarray*}
so applying Theorem~\ref{tsic} with $a=L/\sqrt{t},$ $b=1/\sqrt{t},$ $\alpha=\gamma L/t$ and 
$\beta=\gamma /t$, gives
\begin{eqnarray*}
F_\gamma =
\left(\frac{2\pi}{\sqrt{t}^{\,\gamma(N-1)+1}}\right)^N
A_\gamma(N,t)\,B_\gamma(N,t)\,\prod_{j=0}^{N-1}
\frac{{\Gamma(1+(1+j)\gamma)} }{
{\Gamma(1+\gamma)\Gamma\left(\frac{\gamma}{t}-\gamma j\right)}},
\end{eqnarray*}
where
\begin{eqnarray*}
A_\gamma(N,t)=\lim_{L\rightarrow\infty}
\left(
\frac{\displaystyle{L^{\frac{\gamma L}{t}}\left(\frac{\gamma L}{t}\right)^{\frac{\gamma}{t}-\gamma(N-1)-1}}
}{
\displaystyle{(L+1)^{\frac{\gamma}{t}(L+1)-\gamma(N-1)-1}}
}\right)^N
=\left(\left(\frac{\gamma}{t}\right)^{\frac{\gamma}{t}-\gamma(N-1)-1}e^{-\frac{\gamma}{t}}\right)^N,
\end{eqnarray*}
and, again by standard properties of
the gamma function,
\begin{eqnarray*}
B_\gamma(N,t)=\lim_{L\rightarrow\infty}\prod_{j=0}^{N-1}\left(
\frac{\displaystyle{\Gamma\left(\frac{\gamma}{t}(L+1)-(N+j-1)\gamma-1\right)}
}{
\displaystyle{\Gamma\left(\frac{\gamma L}{t}-\gamma j\right) 
\left(\frac{\gamma L}{t}\right)^{\frac{\gamma}{t}-\gamma (N-1)-1}}
}\right)=1.
\end{eqnarray*}
Thus
$$F_\gamma(N,t)=\left(\frac{{2\pi} }{{t^{\tsh (\gamma(N-1)+1)}}}
 e^{-\gamma/t}
\left(\frac{\gamma}{t}\right)^{\frac{\gamma}{t}-\tsh(\gamma(N-1)+1)}\right)^N
\prod_{j=0}^{N-1}\frac{{\Gamma(1+(1+j)\gamma)} }{ {\Gamma(1+\gamma)
\Gamma\left(\frac{\gamma}{t}-\gamma j\right)}.}$$

\medskip
These evaluations for $F_\gamma$ and $G_\gamma$ immediately yield
\begin{eqnarray}\label{ee9}
W_\gamma(N,t)=\frac{F_\gamma(N,t)}{G_\gamma(N,t)}
=\left(\sqrt{2\pi}\,\, e^{-\gamma/t}
\left(\frac{\gamma}{t}\right)^{\frac{\gamma}{t}-\tsh((N-1)\gamma+1)}\right)^N
\prod_{j=0}^{N-1}\frac{1}{\Gamma\left(\frac{\gamma}{t}-\gamma j\right)}.
\end{eqnarray}

We now exhibit $\log W_\gamma(N,t)$ explicitly as a power series
in $t$, with coefficients that are polynomials in $N$.
Let $N=qK$ and $\gamma=1/q$ where $q$ and $K$ are  positive integers. 
This specialization of $N$ does not affect the recoverability
of $\log W_\gamma(x,t)$ from $\log W_\gamma(N,t)$,
since there is an infinite set of such $K$, so the
polynomiality argument still holds.
With these choices for $N$ and $\gamma$ we have
\begin{eqnarray*}
\prod_{j=0}^{qK-1}{\Gamma\left(\frac{1}{tq}-\frac{j}{q}\right)}=
\frac{\Gamma(\frac{1}{tq}) }{ \Gamma(\frac{1}{tq}-K) }
\prod_{j=1}^{K}\left(\prod_{k=0}^{q-1}
\Gamma\left(\frac{1}{tq}-j+\frac{k}{q}\right)\right),
\end{eqnarray*}
and the product over $k$ can be simplified
by Gauss' multiplication theorem for the gamma function,
whence
\begin{eqnarray}\label{ee10}
\prod_{j=0}^{qK-1}{\Gamma\left(\frac{1}{tq}-\frac{j}{q}\right)}=
\frac{\Gamma(\frac{1}{tq}) }{\Gamma(\frac{1}{tq}-K) }
\left(\frac{\sqrt{2\pi}^{q-1}
}{ q^{\frac{1}{t}-\tsh-q(K+1)}
}
\right)^K
\prod_{j=1}^K\Gamma\left(\frac{1}{t}-qj\right).
\end{eqnarray}
But
$$
\prod_{j=1}^K\Gamma\left(\frac{1}{t}-qj\right)=
\frac{\Gamma(\frac{1}{t})^K t^{qK(K+1)/2}
}{\prod_{j=1}^K\prod_{l=1}^{qj}(1-lt)
},
\quad\mbox{and}\quad
\frac{\Gamma(\frac{1}{tq}) }{ \Gamma(\frac{1}{tq}-K) }=
\frac{1}{(tq)^K}\prod_{j=1}^K(1-tqj).
$$
Combining these results with~(\ref{ee10}) and then with~(\ref{ee9}) gives
\begin{eqnarray}\label{ee11}
W_\frac{1}{q}(qK,t)=
\left(\frac{\sqrt{2\pi t}
}{ \Gamma(\frac{1}{t})\,(et)^\frac{1}{t}
}\right)^K
\frac{{\prod_{l=1}^K\prod_{j=1}^{ql}(1-jt)}
}{ {\prod_{j=1}^K(1-tqj)}
}.
\end{eqnarray}

For the Bernoulli numbers $B_j$ given by~(\ref{berdef}),
recall that
\begin{eqnarray*}
\log\left(
\frac{{\Gamma(\frac{1}{t})\,(et)^\frac{1}{t}}  }{ {\sqrt{2\pi t}}}
\right)=\sum_{k\ge1}\frac{{B_{2k}t^{2k-1}} }{ {2k(2k-1)}}
\end{eqnarray*}
as an asymptotic series in $1/t$, and
\begin{eqnarray*}
\sum_{j=1}^{n}j^k=\frac{1}{k+1}\sum_{r=1}^{k+1}\binom{k+1}{r}B_{k+1-r}
(-1)^{k+1-r}n^r,
\end{eqnarray*}
for $k\geq 0$ and $n\geq 1$.
Then applying these results to~(\ref{ee11}) gives, as an
asymptotic expansion in ${1/t}$,
$$\log W_\frac{1}{q}(qK,t)=
-K\sum_{k\ge1}\frac{{B_{2k}t^{2k-1}}}{2k(2k-1)}
+\sum_{j=1}^K\sum_{\delta\ge1}\frac{{t^{\delta}} }{ \delta}q^{\delta}j^{\delta}
-\sum_{l=1}^K\sum_{j=1}^{ql}\sum_{\delta\ge1}\frac{{t^{\delta}}}{ \delta}
j^{\delta}$$

$$=-K\sum_{k\ge1}\frac{{B_{2k}t^{2k-1}}}{2k(2k-1)}
+\sum_{\delta\ge1}\frac{{t^{\delta}}}{\delta(\delta+1)}
\sum_{r=1}^{\delta+1}\binom{\delta+1}{r}B_{\delta+1-r}(-1)^{\delta+1-r}$$
$$\times\left(q^{\delta} K^r-\frac{q^r}{r+1}
\sum_{m=1}^{r+1}\binom{r+1}{m}B_{r+1-m}(-1)^{r+1-m}
K^m
\right).$$
The result follows by polynomiality, upon
replacing $q$ by ${1/\gamma}$ and $K$ by $N\gamma$, thus
giving
an explicit power series presentation of $\log W_\gamma (x,t)$.
\hfill\qed

\medskip
We are now in a position to give an explicit expression for
the coefficient ${\grgen}^s_g (\gamma )$ given in~(\ref{defH}).
\begin{corollary}\label{corHg}
For g$\geq$1, s$\geq$1,
$${\grgen}^s_g(\gamma )=\left\{
\begin{array}{lcccc} 
\displaystyle{
\frac{(g+s-2)!}{g!}(-1)^s\frac{B_g}{2}
\left(\frac{1}{\gamma^g}-\frac{1}{\gamma}\right)},
&\mbox{if $g$ is even,} \\
\displaystyle{ \frac{{(g+s-2)!(-1)^{s+1}} }{ {(g+1)!}}
\left\{ \frac{{(g+1)B_g} }{{\gamma^g}}+
\sum_{r=0}^{g+1}\binom{ g+1}{r}B_{g+1-r}\frac{B_{r}
}{\gamma^{r}}\right\} }, &\mbox{if $g$ is odd.}
\end{array}
\right.
$$
\end{corollary}
\proof  The ring of asymptotic series in ${1/t}$ is a subring
of the ring of formal power series in $t$. Thus $\log W_{\gamma}$ has
an asymptotic expansion that is the generating series
for $\xi ^s_g (\gamma)$ in ${1/ t}$. Moreover, if a function
has an asymptotic series then it is unique. Thus from
Proposition~\ref{p1} and Theorem~\ref{l3} we obtain,
for $g\geq 1$ and $s\geq 2$,
\begin{eqnarray*}
{\grgen}^s_g (\gamma )&=&s!(-1)^s[x^s t^{g+s-1}]
\frac{1}{\gamma}\log W_\gamma (x,t)\nonumber\\
&=&\frac{s! (-1)^{g+s}}{(g+s)(g+s-1)}
\left\{ \binom{g+s}{s}\frac{B_g}{\gamma ^g}
+\sum_{r=1}^{g+s}
\binom{g+s}{r}B_{g+s-r}\binom{r+1}{s}\frac{B_{r+1-s}}{(r+1)\gamma ^{r+1-s}}
\right\}\nonumber\\
&=&\frac{(g+s-2)!(-1)^{g+s}}{(g+1)!}
\left\{ \frac{(g+1)B_g}{\gamma ^g}
+\sum_{r=1}^{g+s}
\binom{g+1}{r+1-s}B_{g+s-r}\frac{B_{r+1-s}}{\gamma ^{r+1-s}}
\right\}
.\label{Asum}
\end{eqnarray*}
For $g\geq 1$ and $s=1$, we obtain
\begin{eqnarray}
{\grgen}^s_g (\gamma )&=&\frac{(-1)^{g+1}}{g(g+1)}\left\{
B_{g+1}+\frac{(g+1)B_g}{\gamma ^g}
+\sum_{r=1}^{g+1}
\binom{g+1}{r}B_{g+1-r}\frac{B_r}{\gamma ^r}
\right\},\label{Bsum}
\end{eqnarray}
where, for the first term in~(\ref{Bsum}), we have used the
fact that the Bernoulli number $B_j$ is $0$ for $j$ odd and greater
than $1$.

Thus, for all $g\geq 1, s\geq 1$, we have
\begin{eqnarray}
{\grgen}^s_g (\gamma )&=& 
\frac{(g+s-2)!(-1)^{g+s}}{(g+1)!}
\left\{ \frac{(g+1)B_g}{\gamma ^g}
+\sum_{r=0}^{g+s}
\binom{g+1}{r+1-s}B_{g+s-r}\frac{B_{r+1-s}}{\gamma ^{r+1-s}}
\right\}
.\label{Csum}
\end{eqnarray}
(For $s\geq 2$, the summand in (\ref{Csum}) corresponding to $r=0$ is
zero, so (\ref{Csum}) agrees with (\ref{Asum}) in this case.
For $s=1$, the summand in (\ref{Csum}) corresponding to $r=0$ is
 $B_{g+1}$, since $B_0=1$, so (\ref{Csum}) agrees with (\ref{Bsum})
in this case also.) This gives the result immediately for $g$ odd.

For $g$ even, the zero values for Bernoulli numbers of odd subscript
bigger than $1$ mean that for the summation in (\ref{Csum}), only
the terms corresponding
to $r=s$ and $r=g+s-1$ are non-zero, so for $g$ even we obtain
$${\grgen}^s_g (\gamma )=\frac{(g+s-2)!(-1)^s}{(g+1)!}
\left( \frac{(g+1)B_g}{\gamma ^g}
+(g+1)B_g B_1 \frac{1}{\gamma}+(g+1)B_1 B_g\frac{1}{\gamma^g}\right),$$
and the result follows in this case, since $B_1 =-1/2$.
\hfill\qed

\subsection{Proof of Theorem~\ref{main}}
Next we prove Theorem~\ref{main} that
gives ${\Lambda}^s_g$ and ${\Lambda}^{s,O}_{(g+1)/2}$ by
specializing the parameter $\gamma$ in Corollary~\ref{corHg}.
This is straightforward, although it requires additional properties of
the Bernoulli numbers.
\begin{corollary}\label{corchig}
For $g,s\geq1,$
\begin{eqnarray*} 
&\mbox{}&
\begin{array}{ccc}
{\Lambda}^s_g &=&\left\{
\begin{array}{ccccc} 
\displaystyle{\frac{(g+s-2)!}{g!}
(-1)^s (2^{g-1}-1)B_g,}&\mbox{ if $g$ is even,}\\
\displaystyle{\frac{(g+s-2)!(-1)^s}{(g+1)(g-1)!}{B_{g+1}},}
&\mbox{ if $g$ is odd;}
\end{array}
\right.
\end{array}
\\
&\mbox{}&
\begin{array}{ccc}
{\Lambda}^{s,O}_{(g+1)/2} &=&\left\{
\begin{array}{ccccc} 
\displaystyle{0\phantom{\frac{B_g}{g}},}&\mbox{ if $g$ is even,}\\
\displaystyle{ \frac{(g+s-2)!(-1)^s}{(g+1)(g-1)!}B_{g+1}
,}&\mbox{ if $g$ is odd.}
\end{array}
\right.
\end{array}
\end{eqnarray*}
\end{corollary}
\proof From Proposition~\ref{p1} we
have ${\Lambda}^s_g=\grgen^s_g(\tsh)$ and ${\Lambda}^{s,O}_{(g+1)/2}=
{\grgen}^s_g(1)$, so
for $g$ even the result follows immediately from Corollary~\ref{corHg}.

For $g$ odd, we obtain from Corollary~\ref{corHg} that
\begin{equation}\label{Bconv}
\grgen^s_g(\gamma)=\frac{(g+s-2)!(-1)^{s+1}}{(g+1)!}
\left(\frac{(g+1)B_g}{\gamma^g}+
\left[ \frac{t^{g+1}}{(g+1)!}\right] B(t)B(\frac{t}{\gamma})\right),
\end{equation}
where $B(t)$ is the exponential generating series for
Bernoulli numbers, given in~(\ref{berdef}). Now the
following differential equations can be easily verified:
\begin{eqnarray*}
B(t)^2 &=&(1-t)B(t)-t\frac{t}{dt}B(t),\\
B(t)B(2t)&=&(1-t)B(t)-t\frac{t}{dt}B(t)-\frac{t}{2}B(2t).
\end{eqnarray*}
Thus (for $g$ odd), applying the first of these differential equations
to~(\ref{Bconv}) gives
\begin{eqnarray*}
{\grgen}^s_g(1)&=&\frac{(g+s-2)!(-1)^{s+1}}{(g+1)!}
\left\{ (g+1)B_g+\left(
B_{g+1}-(g+1)B_g-(g+1)B_{g+1}\right)\right\} \\
&=&\frac{(g+s-2)!(-1)^s}{(g+1)(g-1)!}
B_{g+1}.
\end{eqnarray*}
Applying the second differential equation gives
\begin{eqnarray*}
{\grgen}^s_g(\tsh )&=&\frac{(g+s-2)!(-1)^{s+1}}{(g+1)!}\left\{
(g+1)B_g 2^g+
\left( B_{g+1}-(g+1)B_g-(g+1)B_{g+1}-(g+1)B_g 2^{g-1}\right)\right\} \\
&=&\frac{(g+s-2)!(-1)^{s+1}}{(g+1)!}\left\{
-g B_{g+1}+(g+1)B_g(2^g-2^{g-1}-1)\right\} ,
\end{eqnarray*}
and the result follows, since $B_g=0$ for $g$ odd except when
 $g=1$, and $2^g-2^{g-1}-1=0$ when $g=1$.
\hfill\qed

\medskip
As discussed at the end of Section~\ref{Stsu}, Corollary~\ref{corchig}
completes the proof of Theorem~\ref{main}.
Note that the natural zero value for ${\Lambda}^{s,O}_{(g+1)/2}$ for  $g$ even
is in agreement with the computed value for $\xi_g^s(1)$  in Corollary~\ref{corchig}.  

\section{A geometric parametrization of the virtual Euler characteristic}\label{Scrvec}

The argument $\gamma$ of the parametrized Euler characteristic
 $\xi_g^s (\gamma )$ was introduced as an artifact for interpolating
between complex curves ($\gamma =1$) and real curves ($\gamma =1/2$),
as  specialized in~(\ref{specpar}). However, we
conjecture that the parameter $\gamma$
itself has geometric significance. In particular, $\xi_g^s (\gamma )$
has been expressed as a polynomial in $1/\gamma$ in Corollary~\ref{corHg},
and we conjecture that the coefficients in this polynomial have a geometric interpretation.

The evidence for this is indirect, but can be lifted from a combinatorial
treatment of maps in Goulden and Jackson~\cite{GJJ}, that specializes
to the series expansion (compare with~(\ref{ee2a})),
\begin{equation}\label{mapgam}
M_{\gamma}(\bfy,x,z) = \sum_{\bfi,j,n}m_{\gamma}(\bfi,j,n)\bfy^\bfi x^jz^n,
\end{equation}
of $M_{\gamma}(\bfy,x,z)$, defined in terms of an integral in~(\ref{ee4}).
Note that $M_{\gamma}(\bfy,x,z)$ is in fact a parametrized map series
because of the specializations in~(\ref{specgam}).
We conjecture that $m_{\gamma}(\bfi,j,n)$ is a polynomial in $1/\gamma$,
with integer coefficients, and that these coefficients have combinatorial
significance.
Moreover,
$$\xi_g^s (\gamma )=s!(-1)^s[x^s t^{g+s-1}]\Psi M_{\gamma},$$
from the proof of Proposition~\ref{p1}, so we deduce that $\xi_g^s (\gamma )$ is
a finite alternating summation of $m_{\gamma}(\bfi,j,n)$'s, and thus we
can lift the
investigation of the significance of $m_{\gamma}(\bfi,j,n)$ to
 $\xi_g^s (\gamma )$ itself. 

The lifting from maps is through properties of symmetric functions,
and Jack functions in particular, which we now summarize (for further
details see Macdonald~\cite{MAC}).

\subsection{Jack symmetric functions}\label{SSJAC}

We say that $\mu=(\mu_1 ,\ldots ,\mu_k )$ is a partition of $n$,
written $\mu\vdash n$ or $\vert\mu\vert =n$, with $k$ parts,
written $l(\mu)=k$, if $\mu_1 \geq\ldots\geq\mu_k \geq 1$,
and $\mu_1 +\cdots +\mu_k =n$. The $\mu_i$'s are the parts of $\mu$.
The monomial symmetric function $m_\mu$ in a countable set of
algebraically independent indeterminates $\bfx =(x_1,\ldots )$,
is the sum of all distinct monomials in the $x_i$'s whose exponents
are the parts of $\mu$, with repetition, in some permuted order.
If $\cP$ is the set of all partitions (including a single, empty,
partition of $0$), then $\{ m_\mu\}_{\mu\in\cP}$ forms a basis
for the ring of symmetric functions. If $p_{\mu}=p_{\mu_1}\cdots p_{\mu_k}$,
where $p_m$ is the $m$th power sum (symmetric function) in $\bfx$,
then $\{ p_{\mu}\} _{\mu\in\cP}$ is also a basis, and hence we can
define an inner product $\langle\phantom{m} ,\phantom{m} \rangle_\alpha$ for an
indeterminate $\alpha$, by the orthogonality
$$\langle p_\lambda,p_\mu\rangle_\alpha
=\frac{\vert\mu\vert!}{\vert{\cal C}_{\mu}\vert}
\alpha^{l(\mu)}\delta_{\lambda,\mu},$$
where ${\cal C}_{\mu}$ is the conjugacy class in the symmetric
group on $\vert\mu\vert$ symbols with disjoint cycle lengths given
by the parts of $\mu$. Then with respect to this inner product, the
Jack symmetric functions $\{ J_{\mu}(\bfx ;\alpha )\}_{\mu\in\cP}$
also form a basis for symmetric functions in $\bfx$ (depending on
the parameter $\alpha$), defined by the conditions
\begin{eqnarray*}
\langle J_\lambda,J_\mu\rangle_\alpha=0\mbox{ for } \lambda\neq\mu;\qquad
[m_\mu]J_\lambda=0\mbox{ unless }\mu\preceq \lambda; \qquad
[x_1\cdots x_n]J_\mu =n!,\mbox{ where } \vert\mu\vert=n,
\end{eqnarray*}
imposing {\em orthogonality, triangularity} and {\em normalization,} respectively,
where $\preceq$ denotes reverse lexicographic order. The Cauchy theorem
for Jack functions is
$$\prod_{i,j\geq 1} (1-x_i y_j )^{-1}=
\sum_{\theta\in\cP}
\frac{J_\theta (\bfx ;\alpha )J_\theta (\bfy ;\alpha )}
{\langle J_\theta ,J_\theta \rangle_\alpha}.$$

Finally, let
$$\langle f(\lambda )\rangle_{\reals ^N} =
\frac{
\displaystyle{\rintn\vert V(\lambda)\vert^{2\gamma}e^{-\frac{\gamma}{2}
p_2}f(\lambda)d\lambda}
}{
\displaystyle{\rintn\vert V(\lambda)\vert^{2\gamma}e^{-\frac{\gamma}{2}
p_2}d\lambda}}.
$$
Then, as a connection between integration and Jack symmetric functions,
it has been previously conjectured in Goulden and Jackson~\cite{GJI} that
\begin{equation}\label{jacsplit}
\langle J_\theta (\lambda ;\alpha)\rangle _{\reals ^N}
= J_\theta (1_N;\alpha )[p_2^m]J_\theta ,
\end{equation}
where $1_N$ is the vector with $N$ $1$'s, and $\theta\vdash 2m$.
This conjecture was recently proved by Okounkov~\cite{OK}.

\subsection{Lifting from the combinatorial conjecture}\label{SSSYM}
We can now give a symmetric function representation for $M_\gamma ({\bf p}(\bfy ),x,z).$

\begin{proposition}\label{intsymp}
Let ${\bf p}(\bfy )=(p_1 (\bfy ),\ldots ).$ Then
\begin{equation}\label{mapjac}
M_\gamma ({\bf p}(\bfy ),x,z)=
\frac{2}{\gamma}z \frac{\partial}{\partial z}
\log\sum_{\theta\in\cP} z^{\frac{\vert\theta\vert}{2}}
\frac{{J_\theta (\bfy ;{1/\gamma})J_\theta (1_N;{1/\gamma})}
}{\langle J_\theta ,J_\theta \rangle_{1/\gamma}}
\left[p_2^{\frac{\vert\theta\vert}{2}}\right]J_\theta .
\end{equation}
\end{proposition}
\proof
If we replace $y_i$ by $p_i (\bfy )$, $i\geq 1$ (this presents no
difficulties since the $p_i (\bfy )$ are algebraically independent
for countable $\bfy$), then from~(\ref{ee4}) and~(\ref{ee5}) we obtain
\begin{eqnarray*}
M_\gamma ({\bf p}(\bfy ),x,z)&=&\frac{2}{\gamma}z \frac{\partial}{\partial z}
\log\left\langle e^{\gamma\sum_{k\geq 1}
\frac{1}{k}\sqrt{z}^k p_k (\bfy )p_k (\lambda )}\right\rangle _{\reals ^N} \\
&=&\frac{2}{\gamma}z \frac{\partial}{\partial z}\log
\left\langle\prod_{i\geq 1}\prod_{j=1}^N
(1-\sqrt{z} y_i \lambda _j)^{-\gamma}\right\rangle _{\reals ^N} \\
&=&\frac{2}{\gamma}z \frac{\partial}{\partial z}\log
\left\langle\sum_{\theta\in\cP}z^{\frac{\vert\theta\vert}{ 2}}
\frac{{J_\theta (\bfy ;{1/\gamma})J_\theta (\lambda ;{1/\gamma})}
}{{\langle J_\theta ,J_\theta \rangle_{1/\gamma}}}\right\rangle _{\reals ^N},
\end{eqnarray*}
from Cauchy's theorem above, and the result follows from the 
integral evaluation~(\ref{jacsplit}).
\hfill\qed

\medskip
Jack function series like that on the right hand side of~(\ref{mapjac}) have been
considered in~\cite{GJJ}, where a combinatorial conjecture is made about
their coefficients. The conjecture specialized to $m_\gamma(\bfi ,j,n)$ is
that $m_{1/(b+1)}(\bfi ,j,n)$ is a polynomial in $b$ with nonnegative
integer coefficients. Note from~(\ref{ee2a}), (\ref{specgam})
and (\ref{mapgam}) that these coefficients sum to $m(\bfi ,j,n)$ (here $b=1$),
and the
constant term is $m(\bfi ,j,n)$ (here $b=0$).
Thus $b$ in this context is a parameter
of nonorientability.

Table 1 of the Appendix gives the values of $m_{1/(b+1)}(\bfi ,j,n)$ for
maps with at most three edges.
For example, the number $4$ in~(\ref{klein}) of  Section 1 is consistent
with the value $m((0,0,0,1),1,2;b)=1+b+3b^2$ in the table. The constant
term $1$ identifies the single orientable rooted map (in the
torus) with $2$ edges,
a single face and a single vertex. The other terms $b$ and $3b^2$
 mean that there are indeed $4$ nonorientable rooted maps (in
the Klein bottle) with $2$ edges,
 $1$ face and $1$ vertex. In this case,
these $4$ maps would be further subdivided into two classes, of size $1$ and $3$,
with different values, $1$ and $2$, respectively, of the unknown measure of
nonorientability recorded by this parameter.

We conclude that $\gamma$, which interpolates between complex curves and
real curves in the context of Euler characteristics, has a separate and
classical existence in terms of the Jack parameter. It should be noted
that, even in the formal study of Jack functions themselves (rather than, for
example, the complicated summation we are considering here), a combinatorial
interpretation
of the Jack parameter has been sought over the last decade (see
Hanlon~\cite{HAN} and Stanley~\cite{ST}), with recent success reported by
Lapointe and Vinet~\cite{VIN} and  Knop and Sahi~\cite{ref21}.

\subsection{Schur symmetric and zonal polynomials}\label{SSZSCH}
We have introduced symmetric functions here indirectly, for technical
reasons. However,  the special cases $\gamma =1/2 ,1$ of the
series $M_\gamma$ give the map series $M,M^O$, respectively, from~(\ref{specgam}).
Moreover, the special cases $\gamma =1/2 ,1$ of Jack functions with
parameter $1/ \gamma$ are zonal polynomials, Schur symmetric functions,
respectively (in the latter case there is also a known scalar introduced).
Applying these specializations
to the Jack function expression~(\ref{mapjac}) therefore yields
expressions with zonal polynomials and Schur
symmetric functions, respectively. These expressions have been
obtained directly by Jackson and Visentin~\cite{JV} (for $\gamma =1$)
and Goulden and Jackson~\cite{GJ2} (for $\gamma =1/2$). (This means
that~(\ref{jacsplit}) and~(\ref{mapjac}) have been proved in the cases
 $\alpha =2,1$ and $\gamma =1/2 ,1$, respectively.) The combinatorial
encodings providing
these direct derivations is quite different from the Wick's lemma
methodology that yielded the integral expressions~(\ref{locint})
for $M$ and~(\ref{orint}) for $M^O$,
in~\cite{GJI} and~\cite{JI}, respectively. We conclude with a sketch
of these combinatorial encodings. For more complete details
see the above references.

A graph embedded in an orientable
surface is rooted, and the two ends of its $n$ edges are labelled from $1$ to
$2n$, with the root end of the
root edge labelled $1$. The map is encoded
by two permutations in $\symgp_{2n}$. In the first, the disjoint
cycles are the cyclic lists of the labels
at the ends of edges incident with a vertex, in clockwise order,
one for each vertex.
In the second, the cycles are transpositions interchanging the
labels at the two ends of a vertex, on for each edge.
The cycles of the product of these permutations
then specify the labels of edges in the boundary of each face of the map, and
the determination of  $M^O$ is therefore reduced to a computation in the centre
of the group algebra of $\symgp_{2n})$,
derivable from the connection coefficients of this subalgebra. The generating
series for these yields the Schur symmetric function expression that
follows  from the summation over $\theta$ in~(\ref{mapjac}) by
setting $\gamma =1$. The logarithm is applied by the standard combinatorial
construction to recover the connected components in this enumeration (since
an arbitrary pair of permutations in the group algebra corresponds
in general to an unordered
collection of labelled maps). The operator $2z{\partial / \partial z}$ is
applied to account for the rooting of the map.

For $M$, a graph embedded in a locally orientable surface has the
four side-end positions of its $n$ edges labelled from $1$ to $4n$,
with the root side and end of the root edge labelled $1$. The map is encoded
by three matchings (pairings) on the $4n$ symbols. In the first, the pairs
consist of the labels that appear at the two ends of the same side of an
edge. In the second, the pairs consist of the labels that appear at the two
sides of the same end of an edge. In the third, the pairs consist of
the labels that appear in a corner of a face.
It can then be shown that the determination of $M$ can therefore be
reduced to a computation in
the double coset algebra of the hyperoctahedral group
embedded in $\symgp_{4n}$ as the centralizer of a fixed fixed point free involution.
The generating series that follows from this yields the zonal polynomial
expression that follows  from the summation over $\theta$ in~(\ref{mapjac}) by
setting $\gamma =1/2$. The logarithm and
application of $4z{\partial / \partial z }$ occur as in the determination
of $M^O$.

\section*{Acknowledgements}
This work was supported by grants individually
to IPG and DMJ from the Natural
Sciences and Engineering Research Council of Canada,
and to JLH from the National Science Foundation (DMS-9401611).


\pagebreak
\appendix
\section{Tables}\label{Atables}
\begin{table}[h]
\begin{eqnarray*}
\begin{array}[t]{lllc|}
n & j & \bfi & m_{1/(b+1)}(\bfi ,j,n) \\ \hline
1 & 1 & (2) & 1 \\
  &   & (0,1) & b \\ \hline
  & 2 & (0,1) & 1 \\  \hline
2 & 1 & (2,1) & 2 \\
  &   & (0,2) & b \\
  &   & (1,0,1) & 4b \\
  &   & (0,0,0,1) & 1+b+3b^2 \\ \hline
  & 2 & (0,2) & 1 \\
  &   & (1,0,1) & 4 \\
  &   & (0,0,0,1) & 5b \\ \hline
  & 3 & (0,0,0,1) & 2 \\ \hline
3 & 1 & (2,2) & 3 \\
  &   & (0,3) & b \\
  &   & (3,0,1) & 2 \\
  &   & (1,1,1) & 12b \\
  &   & (0,0,2) & 1+b+5b^2 \\
  &   & (2,0,0,1) & 9b \\
  &   & (0,1,0,1) & 3+3b+9b^2 \\
  &   & (1,0,0,0,1) & 6+6b+18b^2 \\
  &   & (0,0,0,0,0,1) & 13b+13b^2+15b^3 \\ \hline
\end{array}
\begin{array}[t]{lllc}
n & j & \bfi & m_{1/(b+1)}(\bfi ,j,n) \\ \hline
3 & 2 & (0,3) & 1 \\
  &   & (1,1,1) & 12 \\
  &   & (0,0,2) & 9b \\
  &   & (2,0,0,1) & 9 \\
  &   & (0,1,0,1) & 15b \\
  &   & (1,0,0,0,1) & 30b \\
  &   & (0,0,0,0,0,1) & 10+10b+32b^2 \\ \hline
  & 3 & (0,0,2) & 4 \\
  &   & (0,1,0,1) & 6 \\
  &   & (1,0,0,0,1) & 12 \\
  &   & (0,0,0,0,0,1) & 22b \\ \hline
  & 4 & (0,0,0,0,0,1) & 5 \\ \hline
\end{array}
\end{eqnarray*}
\caption{The refined map
numbers $m_{1/(b+1)}(\bfi ,j,n)$, for $n\leq 3$ edges.}
\end{table}

\end{document}